\documentclass[review]{elsarticle}

\journal{Automatica}

\usepackage{psfrag}
\usepackage{amsthm}
\usepackage{amsmath}
\usepackage{amssymb}
\usepackage{subfigure}

\newtheorem{Theorem}{Theorem}
\newtheorem{Definition}{Definition}
\newtheorem{Lemma}{Lemma}

\newcommand{\sN}{\mathcal{N}}

\newcommand{\abs}[1]{\left\vert #1 \right\vert}
\newcommand{\norm}[1]{\left\Vert #1 \right\Vert}

\newcommand{\R}{{\mathbb R}}  

\begin{document}

\begin{frontmatter}

\title{On noise-induced synchronization and consensus}

\author[GRmainaddress]{Giovanni Russo\corref{mycorrespondingauthor}}
\cortext[mycorrespondingauthor]{Corresponding author}
\ead{grusso@ie.ibm.com}
\author[GRmainaddress,RSsecondaryaddress]{Robert Shorten}
\address[GRmainaddress]{IBM Research Ireland, Optimization, Control and Decision Science Group}
\address[RSsecondaryaddress]{University College Dublin, School of Electrical and Electronic Engineering}

\begin{abstract}
In this paper, we present new results for the synchronization and consensus of networks described by Ito stochastic differential equations. From the methodological viewpoint, our results are based on the use of stochastic Lyapunov functions. This approach allowed us to consider networks where nodes dynamics can be nonlinear and non-autonomous and where noise is not just additive but rather its diffusion can be nonlinear and depend on the network state. We first present a sufficient condition on the coupling strength and topology ensuring that a network synchronizes (fulfills consensus) {\em despite} noise. Then, we show that noise can be useful, and present a result showing how to {\em design} noise so that it induces synchronization/consensus. Motivated by our current research in Smart Cities and Internet of Things, we also illustrate the effectiveness of our approach by showing how our results can be used to analyze/control the onset of synchronization in noisy networks and to study collective decision processes.
\end{abstract}


\end{frontmatter}


\section{Introduction}

Network control is of utmost importance in many application fields, ranging from computer science to power networks, internet of things and systems biology \cite{Cor_Kat_Mot_13}, \cite{Liu_Bar_Slo_11}. In all such fields, the problem of steering the dynamics of network agents towards a coordinated behavior is recognized as a fundamental network control problem, \cite{Che_13}. Synchronization and consensus are two examples where all the agents of the network need to coordinate their actions in a decentralized way so that they all converge to the same behavior, \cite{Olf_Mur_04}, \cite{Su_Hua_12}, \cite{Dor_Bul_14}. Consensus and synchronization are also at the core of distributed optimization algorithms, see e.g. \cite{Dro_Kaw_Ege_14}, \cite{Stu_Sho_12} and references therein, where a set of agents {\em collaborates} to minimize a global cost function by only using some locally available information.

Over the past few years, several sufficient conditions have been devised ensuring, under different technical assumptions, synchronization/consensus of complex networks. However, an assumption that is often made in Literature is that the network of interest is {\em noise-free}. This assumption is not realistic for most real world applications of synchronization and consensus, where noise plays a key role in destroying or generating those coordinated behaviors. Networks of power generators communicating over transmission lines \cite{Dor_Bul_14}, neural networks \cite{Wil_Moh_14a}, distributed estimation \cite{Zha_Dua_Wen_Che_16} and cognitive processes \cite{Pet_Sri_Tay_Sur_Eck_Bul_15} are all examples of applications where noise cannot be neglected. 

Motivated by this, the problem of controlling synchronization/consensus in networks affected by noise has recently attracted many researchers. For example, in \cite{Tan_Li_15} the consensus problem is studied when the information exchanged among the network nodes is corrupted, while in \cite{Sri_Leo_14} coupled drift-diffusion models are used to investigate the dynamics of collective decision making processes. In all these papers, it is assumed that the noise affecting each node is additive with constant diffusion rate and that the dynamics at the network nodes are linear systems or integrators. Recently, an approach to study synchronization in networks of nonlinear nodes subject to additive noise has been presented in \cite{Wel_Kat_Mot_15}. However, the key idea behind the methodology presented in this paper is to recast synchronization as a stochastic optimization problem and rely on numerical methods to solve it.

In this paper, we present new algebraic sufficient conditions ensuring that synchronization/consensus is attained for a given network. We will consider a wide class of networks relevant to applications where nodes' dynamics are nonlinear and where the noise affecting nodes has a diffusion rate that can also be nonlinear and state dependent. By using Lyapunov techniques \cite{Mao_97}, we will first present a sufficient condition ensuring that a network synchronizes (or achieves consensus) in the presence of noise. Our condition will explicitly link together network topology, coupling strength, node and diffusion dynamics and, intuitively, it will imply that network will achieve synchronization/consensus if the coupling strength and the network topology are well {\em blended} together so as to overcome a {\em threshold} generated by noise.

Then, we will show that, if properly {\em designed}, noise can be also turned into a {\em distributed control input} to induce synchronization/consensus. Specifically, we will show that a given network of interest can be driven towards synchronization/consensus by designing an additional layer, where noise is {\em injected}. In essence, this new layer, which superimposes to the original network topology, controls the network by properly canalizing noise. Interestingly, the new {\em noise-propagation layer} can have a topology different from the original network topology and, in principle, nodes might  even be disconnected at the network level but connected through this newly designed layer. Recently, multi-layer network control has been investigated in \cite{Bur_diB_15} but in this case noise was seen as a disturbance to be rejected through the development of a distributed PI controller. Our viewpoint is instead fundamentally different as we aim to {\em embrace} noise and make it useful for network control. To the best of our knowledge, the idea of designing an additional layer where noise can propagate to control a network is new in this Literature. In developing this idea we have been inspired by our current research in Smart Cities and Internet of Things (IoT). For example, for such applications, it often the case where, due to physical/environmental constraints, the coupling strength between nodes needs to be low, and thus such networks, generally speaking, are unable to achieve consensus/synchronization. In this case, our results can be used to properly inject some noise that, combined with the existing couplings, induces synchronization. Another motivation, particularly relevant of applications of control to {\em IoT}, is that typically network nodes need to preserve some level of {\em privacy}. This is, for example, the case in vehicles platooning, where a common speed needs to be found among a network of vehicles but where each driver is not willing to share his/her exact position, speed and utility function with the other drivers. In this case, with our results we can use noise in order to {\em obfuscate} private data so that no node knows the exact state of its neighbors and, at the same time, network coordination is achieved.

The paper is organized as follows. We start in Section \ref{sec:math_prel} by introducing the notation and mathematical preliminaries that are useful for the paper, while in Section \ref{sec:probl_stat} we will formalize the problem statement. Then, in Section \ref{sec:com_sto_syn} we introduce a sufficient condition to ensure complete stochastic synchronization of a network affected by noise, while in Section \ref{sec:noise_synch} we turn our attention to network control and provide a result that allows to use noise to achieve synchronization. The effectiveness of our approach is demonstrated in both Section \ref{sec:com_sto_syn} and Section \ref{sec:noise_synch} where our results are used to analyze/control synchronization in a network of Fitzhugh-Nagumo oscillators. Finally, in Section \ref{sec:applications} we also show that the approach proposed in this paper can be applied to study a collective decision making (or {\em cognitive}) processes.

\section{Mathematical Preliminaries}\label{sec:math_prel}

\subsection{Notation}
In this paper, we denote by $I_n$ the $n\times n$ identity matrix and by $1_{n\times m}$ the $n \times m$ matrix having all of its elements equal to $1$. The vector/matrix Frobenius norm will be denoted by $\norm{\cdot}_F$ and the vector/matrix Euclidean norm will be denoted by $\abs{\cdot}$. The trace of a square matrix, say $A$, will be denoted by $tr\left\{A\right\}$.

\subsection{Stochastic differential equations}
Consider an $n$-dimensional stochastic differential equation of the form
\begin{equation}\label{eqn:ito_gen}
dx = f(t,x)dt + g(t,x)dB,
\end{equation}
where: (i) $x \in \R^n$ is the state variable; (ii) $f:\R^+\times\R^n\rightarrow\R^n$ belongs to $\mathcal{C}^2$; (iii) $g:\R^+\times\R^n\rightarrow\R^{n\times d}$ belongs to $\mathcal{C}$; (iv) $B= [B_1,\ldots,B_d]^T$ is a $d$-dimensional Brownian motion. Throughout this paper we will assume that both $f$ and $g$ obey the local Lipschitz condition and the linear growth condition, see e.g. \cite{Oks_07}. This implies that for any given initial condition $x(t_0) = x_0$, $t\ge 0$, equation (\ref{eqn:ito_gen}) has a unique global solution. We will also assume that $f(t,0) = g(t,0) = 0$ and the solution $x=0$ will be said the {\em trivial solution} of (\ref{eqn:ito_gen}).

Following \cite{Kar_93}, \cite{Roh_76}, we say that a sequence of stochastic variables, $\left\{V_1,V_2,\ldots\right\}$ converges almost surely (a.s.) to the stochastic variable $V$ if
$$
\mathbb{P}\left(\left\{w: \lim_{n\rightarrow +\infty}V_n(w)=V(w)\right\}\right) = 1.
$$
That is, the sequence converges to $V$ with probability $1$. We are now ready to give the following definition which characterizes stability of the trivial solution, see \cite{Mao_97}. 
\begin{Definition}
The trivial solution of (\ref{eqn:ito_gen}) is said to be almost surely exponentially stable if for all $x\in\R^n$, $\lim_{t\rightarrow +\infty}\sup\frac{1}{t}\log\left(\abs{x(t)}\right) <0$, $a.s.$.
\end{Definition}

Let $V(t,x):\R^+\times\R^n\rightarrow\R^+$, $V(t,x) \in \mathcal{C}^{1\times2}$, i.e. $V(t,x)$ is twice differentiable with respect to $x$ and differentiable with respect to $t$. By the Ito formula we have:
$$
dV(t,x) = LV(t,x)dt +V_x(t,x)g(t,x)dB,
$$
where: (i) $LV(t,x) = V_t(t,x) + V_x(t,x)f(t,x) + \frac{1}{2}tr\left\{ g(t,x)^TV_{xx}g(t,x)(t,x)\right\}$;(ii) $V_x = \left[V_{x_1},\ldots,V_{x_n}\right]$; (iii) $V_{xx}$ is the $n\times n$ dimensional matrix having as element $ij$ $V_{x_ix_j}$ (where $V_{x_i} := \partial V(t,x)/\partial x_i$ and $V_{x_ix_j} := \partial^2 V(t,x)/\partial x_j \partial x_i$). The following result from \cite{Mao_97} provides a sufficient condition for the trivial solution of (\ref{eqn:ito_gen}) to be almost surely exponentially stable.
\begin{Theorem}\label{thm:stability_ito}
Assume that there exists a non-negative function $V(t,x) \in \mathcal{C}^{1\times2}$ and constants $p>0$, $c_1>0$, $c_2\in\R$, $c_3\ge 0$, such that $\forall x \ne 0$ and $\forall t \in \R^+$: ({\bf H1}) $c_1\abs{x}^p \le V(t,x)^p$; ({\bf H2}) $LV(t,x) \le c_2 V(t,x)$; ({\bf H3}) $\abs{V_x(t,x)g(t,x)}^2 \ge c_3 V(t,x)^2$.
Then:
$ \lim_{t\rightarrow +\infty}\sup\frac{1}{t}\log\left(\abs{x(t)}\right) \le -\frac{c_3-2c_2}{p}$, $a.s.$. In particular, if $c_3>2c_2$, then the trivial solution of (\ref{eqn:ito_gen}) is almost surely exponentially stable. 
\end{Theorem}
\subsection{Complex networks}
Throughout this paper, we will consider systems interacting over some graph, $\mathcal{G} = (\mathcal{V},\mathcal{E})$, where $\mathcal{V}$ is the set of vertices (or nodes in what follows) and $\mathcal{E}$ is the set of edges. We assume all the graphs in this paper are undirected and denote the edge between node $i$ and node $j$ as $(i,j)$. We will denote by $\sN_i$ the set of neighbors of node $i$, i.e. $\sN_i :=\left\{j:(i,j)\in\mathcal{E}\right\}$. Let $N$ be the number of nodes in the network. Then (see e.g. \cite{God_Roy_01}) the Laplacian matrix associated to the graph, $L$, is the $N\times N$ symmetric matrix defined as $L = \Delta - A$, where: (i) $A$ is the adjacency matrix of $\mathcal{G}$; (ii) $\Delta$ is the graph degree matrix. The following Lemma from \cite{Hor_Joh_99} will be used in this paper.

\begin{Lemma}\label{lem:laplacian}
Denote by $L$ the Laplacian matrix of an undirected network. The following properties hold: (i) $L$ has a simple zero eigenvalue and all the other eigenvalues are positive if and only if the network is connected; (ii) the eigenvector associated to the zero eigenvalue is $1_N$, i.e. the $N$-dimensional vector having all of its elements equal to $1$; (iii) the smallest nonzero eigenvalue, $\lambda_2$, satisfies $\lambda_2 = \min_{v^T 1_N=0, v \ne 0}\frac{v^TLv}{v^Tv}$.
\end{Lemma}
In the rest of this paper, $\lambda_2$ will be termed as the graph algebraic connectivity. 
\section{Problem statement}\label{sec:probl_stat}
Throughout this paper, we will consider stochastic networks described by the following stochastic differential equation:
\begin{equation}\label{eqn:net_more_general}
dx_i = \left[f(t,x_i) + \sigma\sum_{j\in\sN_i}\left(x_j-x_i\right)\right]dt + g_i(t,X)db_i, \ \ \ i = 1,\ldots,N,
\end{equation}
where: (i) $x_i \in \R^n$; (ii) $f(t,x_i):\R^+\times\R^n\rightarrow\R^n$ is the smooth {\em nominal} nodes dynamics; (iii) $\sigma$ is the coupling strength; (iv) $g_i(t,X): \R^+\times\R^{nN}\rightarrow\R^{n\times d}$ is the smooth {\em diffusion} matrix describing how noise affects node $i$; (v) $b_i(t) = \left[b_{i,1},\ldots,b_{i,d}\right]^T$ is the $d$-dimensional Brownian process describing the noise acting on the $i$-th network node. Network dynamics (\ref{eqn:net_more_general}) can be written in compact form as follows:
\begin{equation}\label{eqn:network}
dX = \left[F(t,X) - \sigma (L\otimes I_n)X\right]dt + G(t,X)dB,
\end{equation}
with: (i) $X = [x_1^T,\ldots,x_N^T]^T$; (ii) $F(t,X) = [f(t,x_1)^T,\ldots,f(t,x_N)^T]^T$; (iii) $B = \left[b_1^T,\ldots,b_N^T\right]^T$; (iv) $G(t,X)$ is the $Nn\times Nd$ block diagonal matrix having  $G_{ii}(t,X) = g_i(t,X)$. Please note that in the case where $db_i = db_j \in \R$, $\forall i,j = 1,\ldots,N$, then the term $G(t,X)dB$ can be written as $G(t,X)db$, with $G(t,X)$ being the $Nn$-dimensional column vector having $G_i(t,X) = g_{i}(t,X)$.

The goal of this paper is to address the so-called {\em synchronization problem}.  This is formalized with the following definition.
\begin{Definition}
Let $\dot s_n(t)= \frac{1}{N}\sum_{i=1}^Nf(t,x_i)$. We will say that network (\ref{eqn:net_more_general}) achieves complete stochastic synchronization if
$\lim_{t\rightarrow +\infty}\sup\frac{1}{t}\log\left(\abs{x_i(t)-s_n(t)}\right) <0$, $a.s.$, $\forall i=1,\ldots,N.$.
\end{Definition}
Note that in the case where nodes' dynamics are integrator dynamics, then the definition of complete stochastic synchronization simply becomes a definition for consensus.
\section{A sufficient condition for complete stochastic synchronization}\label{sec:com_sto_syn}

The following result provides a sufficient condition ensuring complete stochastic synchronization of network (\ref{eqn:net_more_general}).  
\begin{Theorem}\label{thm:network_nodes_general}
Let $S = 1_N\otimes s_n(t)$ and assume that for network (\ref{eqn:net_more_general}) the following conditions are fulfilled:
\begin{enumerate}
\item there exists some constant, say $K_f$, such that $(X-Y)^T\left[F(t,X)-F(t,Y)\right]\le K_f(X-Y)^T(X-Y)$, $\forall X,Y \in \R^{Nn}$;
\item $G(t,S)=0$, $\forall t$;
\item there exists some constant, say $K_g$, such that $\norm{G(t,X)-G(t,S)}_F\le K_g\norm{X-S}_F$, $\forall X \in \R^{Nn}$;
\item there exists some constant, say $\bar K_g$, such that  $\abs{(X-S)^T\left(G(t,X)-G(t,S)\right)}^2\ge \bar K_g^2\abs{X-S}^4$, $\forall X \in \R^{Nn}$;
\item $\sigma\lambda_2 > K_f + \frac{K_g^2-2\bar K_g^2}{2}$.
\end{enumerate}
Then, (\ref{eqn:net_more_general}) achieves complete stochastic synchronization.
\end{Theorem}
\proof
Let $e = X-S$. Then, the error dynamics can then be written as:
\begin{equation}\label{eqn:error_3}
de = \left[\tilde F(t,e)\right]dt + \left[\tilde G(t,e)\right]dB,
\end{equation}
where:
\begin{itemize}
\item $\tilde F(t,e) =F(t,e+S) - \sigma (L\otimes I_n)(e+S) - \frac{1}{N}\left(1_{N\times N}\otimes I_n\right)F(t,e+S)$;
\item $\tilde G(t,e) = G(t,e+S)$. 
\end{itemize}
Note that $e=0$ is the trivial solution for (\ref{eqn:error_3}). In fact:
$F(t,S) - \frac{1}{N}1_{N\times N}F(t,S) = 0$ and, by Hypothesis $2$,
$G(t,S) = 0$.
Let $V(t,e) = V(e) = \frac{1}{2}e^Te$, then following Theorem \ref{thm:stability_ito}, in order to prove our result we need to show that there exists $c_2\in\R$, $c_3\ge 0$, such that $c_3 > 2c_2$. In order to prove this, we will now estimate $LV(e)$ and $\abs{V_e(e)\tilde G(t,e)}^2$.

{\bf Estimate of $LV(e)$.} In order to compute this term, first note that $V_t(e) = 0$. Let's now compute the term $V_e(e)\tilde F(t,e)$. We have:
$$
V_e(e)\tilde F(t,e) = e^T\left[F(t,e+S) - \sigma (L\otimes I_n)(e+S) - \frac{1}{N}\left(1_{N\times N}\otimes I_n\right)F(t,e+S) \right],
$$
and, by adding and subtracting $F(t,S)$, we get
$$
\begin{array}{*{20}l}
V_e(e)\tilde F(t,e) & = & e^T\left[F(t,e+S)-F(t,S) + F(t,S) - \sigma (L\otimes I_n)(e+S) + \right.\\
& + & \left.- \frac{1}{N}\left(1_{N\times N}\otimes I_n\right)F(t,e+S) \right].
\end{array}
$$
On the other hand, note that $
e^T\left[F(t,S) - \frac{1}{N}\left(1_{N\times N}\otimes I_n\right)F(t,e+S) \right] = 0$, while: $- \sigma (L\otimes I_n)(e+S) = -\sigma (L\otimes I_n) e$ and therefore 
\begin{equation}\label{eqn:diffusion}
V_e(e)\tilde F(t,e) = e^T\left[F(t,e+S)-F(t,S)\right] - \sigma e^T\left(L\otimes I_n\right)e. 
\end{equation}
Now:
$$
V_e(e)\tilde F(t,e) \le e^T\left[F(t,e+S)-F(t,S)\right] - \sigma min_{e\ne0}\left\{e^T(L\otimes I_n)e\right\},
$$
and since $e^T 1_{Nn} =0$, then by means of Lemma \ref{lem:laplacian} we have $\min_{e \ne 0}\left[e^T  \left(L \otimes I_n\right)e\right] = \lambda_2 e^T\left(I_N\otimes I_n\right)e = \lambda_2 e^Te$. 

Finally, by Hypothesis $1$ we have $e^T\left[F(t,e+S)-F(t,S)\right] \le K_fe^Te$, and thus:
\begin{equation}\label{eqn:LV_3}
V_e(e)\tilde F(t,e) \le \left(K_f-\sigma\lambda_2\right)e^Te = 2\left(K_f-\sigma\lambda_2\right)V(e).
\end{equation}
The next step to estimate $LV(e)$ is that of computing $\frac{1}{2}tr\left\{ \tilde G(t,e)^TV_{ee}\tilde G(t,e)\right\}$. Now, since $V_{ee}(e) = I_{Nn}$, such a term simply becomes:
$\frac{1}{2}tr\left\{ \tilde G(t,e)^T\tilde G(t,e)\right\}$.
Also, recall that for any matrix, say $A$, we have $\norm{A}_F^2 = tr\left\{ A^TA\right\}$. Thus:
$$
\left(tr\left\{ \tilde G(t,e)^TV_{ee}\tilde G(t,e)\right\}\right)^{1/2} =\norm{\tilde G(t,e)}_F = \norm{G(t,e+S)}_F.
$$
Therefore, by means of Hypothesis $3$ we have:
$$
\left(tr\left\{ \tilde G(t,e)^TV_{ee}\tilde G(t,e)\right\}\right)^{1/2} =\norm{G(t,e+S)-G(t,S)}_F\le K_g \norm{e}_F = K_g \abs{e},
$$
where the last inequality follows from the fact that for any vector, say $a$, $\norm{a}_F = \abs{a}$. We finally have:
\begin{equation}\label{eqn:VXX_3}
\frac{1}{2}tr\left\{ \tilde G(t,e)^TV_{ee}\tilde G(t,e)\right\} \le \frac{1}{2}K_g^2\abs{e}^2 =\frac{1}{2} K_g^2e^Te = K_g^2V(e).
\end{equation}
Combining (\ref{eqn:VXX_3}) and (\ref{eqn:LV_3}) we get:
$$
LV(e) \le \left(2K_f+K_g^2-2\sigma\lambda_2\right)V(e) := c_2 V(e).
$$
Now, in order to complete the proof we need to find a lower bound for $\abs{V_e(e)\tilde G(t,e)}^2$.

{\bf Estimate of $\abs{V_e(e)\tilde G(t,e)}^2$.} In order to estimate this term, first note that
$$
\abs{V_e(e)\tilde G(t,e)} = \abs{e^T\left(G(t,e+S) - G(t,S)\right)}.
$$
Thus, by Hypothesis $4$, we have
\begin{equation}\label{eqn:sec_term}
\abs{V_e(e)\tilde G(t,e)}^2 \ge \bar K_g^2 \abs{e}^4 = \bar K_g \left(e^Te\right)^2 = 4\bar K_g^2V(e)^2:=c_3V(e)^2.
\end{equation}

We can then conclude the proof by noticing that, by Hypothesis $5$
$$
4\bar K_g^2 > 2\left(2K_f+K_g^2-2\sigma\lambda_2\right).
$$
Therefore, by means of Theorem \ref{thm:stability_ito}, 
$\lim_{t\rightarrow +\infty}\sup\frac{1}{t}\log\left(\abs{e(t)}\right) <0$, $a.s.$, thus proving the result.
\endproof

\subsection*{Example: synchronization of Fitzhugh-Nagumo oscillators}

In order to illustrate the key features of Theorem \ref{thm:network_nodes_general}, we consider a network of Fitzhugh-Nagumo (FN, see \cite{Fit_55}) oscillators. Specifically, we will consider the case where network nodes do not communicate directly but they rather communicate by means of a shared environmental variable which is affected by some noise. In neuronal contexts, this mechanism is known as \emph{local field potential} and it is believed to play an important role in the synchronization of groups of neurons, \cite{Ana_Mon_Bar_Buz_Koc_10}, \cite{Bre_Lai_11}, \cite{Rus_Slo_10}.

\begin{figure}[tbh]
\begin{center}
\centering
\psfrag{x1}[c]{{$\gamma$}}
\psfrag{x2}[c]{{$t$}}
\psfrag{x3}[c]{{$t$}}
\psfrag{y1}[c]{{$R$}}
\psfrag{y2}[c]{{$v_i$'s}}
\psfrag{y3}[c]{{$v_i$'s}}
  \includegraphics[width=8cm]{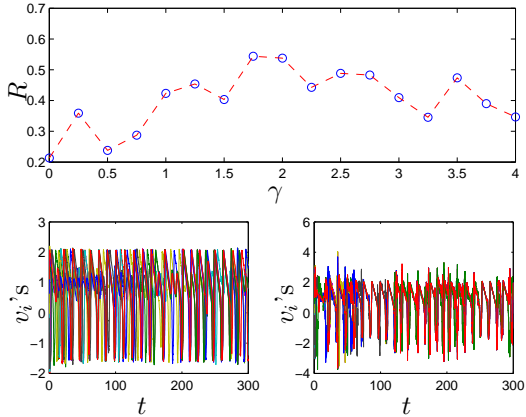}
  \caption{Simulation of (\ref{eqn:FN}) for $N=10$ and the following set of parameters for the nodes' dynamics: $a = 0.7$, $b = 0.4$, $c=2.8$, $u(t) = 0$ (this set of parameters gives $K_f = 3.2$). In the top panel, the behavior of $R$ is shown as a function of the noise intensity, $\gamma$. Network synchronization is not attained and this is confirmed by the bottom panels, where the time behaviors of the network state variables are shown for $\gamma = 0.25$ (left-bottom panel) and $\gamma = 3.75$ (right-bottom panel). The simulations were run using Matlab/Simulink, with solver method {\it ode23s} and relative tolerance {\it 1e-5}.}
  \label{fig:order_param}
  \end{center}
\end{figure}

The following mathematical model has been adapted from \cite{Rus_Slo_10} to describe a network of FN oscillators coupled through a noisy field potential:
\begin{equation}\label{eqn:FN}
\begin{array}{*{20}l}
d v_i = \left[c\left(v_i+w_i -\frac{1}{3}v_i^3 +u(t)\right)\right]dt + \left[\gamma\left(\frac{1}{N}\sum_{j=1}^N(v_j) - v_i\right)\right]db\\
d w_i = \left[-\frac{1}{c}\left(v_i-a+bw_i\right)\right]dt,\\
\end{array}
\end{equation}
where: (i) $v_i$ and $w_i$ are the membrane potential and the recovery variable for the $i$-th FN oscillator ($i=1,\ldots,N$); (ii) $u(t)$ is the smooth magnitude of an external stimulus taken as input by all the FN oscillators in the network; (iii) the term $\left[\sum_{j=1}^N(v_j) - v_i\right]$ describes the interaction between each of the nodes and the shared environment; (iv) $db$ is the standard scalar Brownian motion acting on the environment; (v) $\gamma$ is the noise intensity. In terms of the notation of Theorem \ref{thm:network_nodes_general} we have (please see the Appendix for the computations): (i) $\sigma = \lambda_2 = 0$; (ii) $K_f = \frac{-b+c^2 +\sqrt{1+b^2+2(b-1)c^2+2c^4}}{2c}$; (iii) $K_g = \gamma N$ and $\bar K_g = 0$. Therefore, following Theorem \ref{thm:network_nodes_general}, it follows that network (\ref{eqn:FN}) will synchronize if $K_f + (\gamma N)^2/2 <0$. In Figure \ref{fig:order_param} we simulated network (\ref{eqn:FN}) with $N = 10$ and with $\gamma$ ranging between $0$ and $4$. In order to characterize quantitatively the {\em level} of synchronization of the network we used the {\em order parameter} 
$R := (\langle M^2\rangle - \langle M\rangle^2)/(\overline{\langle v_i^2\rangle - \langle v_i\rangle^2})$, defined follwing \cite{Gar_Elo_Str_04} where: (i) $M(t) := 1/N\sum_{i=1}^Nv_i$; (ii) $\langle \cdot \rangle$ denotes the time average; (ii) $\bar{\cdot}$ denotes the average over the network nodes. In Figure \ref{fig:order_param} (top panel) the order parameter for network (\ref{eqn:FN}) is plotted as a function of the noise intensity $\gamma$. Such a panel shows that network synchronization is not attained. In the bottom panels of Figure \ref{fig:order_param} the time behaviors of $v_i$'s are instead shown for $\gamma = 0.25$ and $\gamma = 3.75$. We will revisit field potentials in the example of Section \ref{sec:noise_synch} where we will show that, if properly {\em designed}, noise diffusion can induce synchronization.

\section{Remarks}

\begin{itemize}
\item Hypothesis $1$ is sometimes known in the literature as {\em QUAD}. As shown in \cite{deL_diB_Rus_11}, this condition can be linked to Lipschitz and contraction conditions of the vector field (in this latter case, this would imply $K_f < 0$);
\item Hypothesis $2$ implies that the noise disappears when the synchronous state is reached. An important class of functions satisfying this hypotheses are {\em Laplacian-like} functions. Such functions will be used in Section \ref{sec:noise_synch} to synchronize a network through noise;
\item Hypothesis $3$ implies that the noise diffusion is bounded while Hypothesis $4$ allows the noise to be {\em persistent};
\item Hypothesis $5$ links together a structural property of the network graph (i.e. the algebraic connectivity) and the strength of the coupling between nodes. A sufficient condition for synchronization is for $\sigma \lambda_2$ to exceed the threshold $\tilde K := K_f + \frac{K_g^2-2\bar K_g^2}{2}$ which depends on the dynamics of the node (through $K_f$) and on noise (through $K_g$ and $\bar K_g$);
\item Assume that $2\bar K_g^2 > K_g^2$. In this case, the effect of noise diffusing through the network is that of lowering $\tilde K$. That is, noise plays in favor of synchronization. Motivated by this, in Section \ref{sec:noise_synch}, we will present a methodology that will allow to design noise so that a given network of interest synchronizes.
\end{itemize}

\section{Inducing synchronization through noise}\label{sec:noise_synch}
\begin{figure}[tbh]
\begin{center}
  \includegraphics[width=6cm]{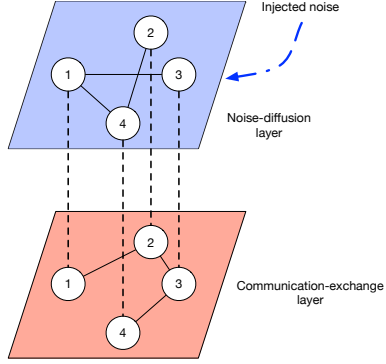}
  \caption{Controlling synchronization through noise. Essentially noise is injected in the top layer and propagates therein. This noise-diffusion layer communicates with the actual network where nodes exchange information (state variables, $x_i$'s). If properly designed, the noise-diffusion layer can enforce network synchronization. Note that that the network topology of the two layers are in general different.}
  \label{fig:layer}
  \end{center}
\end{figure}

In this Section we turn our attention to the problem of {\em designing} noise in order to synchronize a given network of interest. The set-up we have in mind is outlined in Figure \ref{fig:layer}. The key idea of our approach is to design a new {\em network layer} where noise in injected and can propagate. It is then this new layer that induces synchronization. Formally, the resulting network topology is described by a multi-graph \cite{Bur_diB_15}. Here, the multi-graph consists of two layers: (i) the {\em communication} graph topology, where nodes information is exchanged (i.e. nodes exchange their state variables); (ii) the {\em noise-diffusion} layer which has the goal of enforcing network synchronization.

Our set-up is motivated by our current research in Smart Cities and IoT. For example, in networks affected by environmental constraints (i.e. low coupling strengths unable to synchronize the network) the additional noise-diffusion layer can be used to properly inject low intensity noise that, combined with the existing couplings induces synchronization. Another motivation particularly relevant for applications of control to IoT is that network nodes need to {\em preserve their privacy}. In this case, our idea is to use the noise-diffusion layer to {\em blend} private data (nodes' state variables) with noise so that synchronization is still attained while no agent knows the exact state of its neighbors. With respect to this, note that the topologies of the two network layers are in general different and, in principle, nodes might be disconnected at the communication layer and interact with each other only over the noise-diffusion layer.

Formally, the set-up of Figure \ref{fig:layer} corresponds to the following set of stochastic Ito differential equations:
\begin{equation}\label{eqn:net_control}
dx_i = \left[f(t,x_i) + \sigma\sum_{j\in\sN_i}\left(x_j-x_i\right)\right]dt + \sigma^\ast \sum_{j\in\sN_i^\ast}\left(x_j-x_i\right) db, 
\end{equation}
$i= 1,\ldots,N$, which in compact form leads to
\begin{equation}\label{eqn:net_control_compact}
dX = \left[F(t,X) - \sigma (L\otimes I_n)X\right]dt -\sigma^\ast(L^\ast\otimes I_n)Xdb.
\end{equation}
Note that all the {\em structural} quantities of the noise diffusion layer are denoted by the superscript $\ast$. We are now ready to give the following result.
\begin{Theorem}\label{thm:control}
Let:
\begin{itemize}
\item $\lambda_N^\ast$ be the largest eigenvalue of the graph Laplacian of the noise diffusion layer;
\item $\lambda_2^\ast$ be the algebraic connectivity of noise diffusion layer;
\end{itemize}  
and assume that:
\begin{enumerate}
\item there exists some constant, say $K_f$, such that $(X-Y)^T\left(F(t,X)-f(t,Y)\right)\le K_f(X-Y)^T(X-Y)$, $\forall x,y \in \R^{Nn}$;
\item $(\sigma^\ast)^2\left((\lambda_2^\ast)^2-\frac{(\lambda_N^\ast)^2}{2}\right)>K_f-\sigma\lambda_2$.
\end{enumerate}
Then, (\ref{eqn:net_control}) achieves complete stochastic synchronization.
\end{Theorem}
\proof
Following steps similar to those used to prove Theorem \ref{thm:network_nodes_general}, we define the quantity $e = X-S$ and we have:
\begin{equation}\label{eqn:error_noise}
de = \left[\tilde F(t,e)\right]dt + \left[\tilde G(t,e)\right]db,
\end{equation}
where:
\begin{itemize}
\item $\tilde F(t,e) =F(t,e+S) - \sigma (L\otimes I_n)(e+S) - \frac{1}{N}\left(1_{N\times N}\otimes I_n\right)F(t,e+S)$;
\item $\tilde G(e) = -\sigma^\ast (L^\ast\otimes I_n)(e+S) = -\sigma^\ast (L^\ast\otimes I_n)e$.
\end{itemize}
We can again use Theorem \ref{thm:stability_ito} to prove our result with $V(t,e) = V(e) = \frac{1}{2}e^Te$.

{\bf Estimate of $LV(e)$.} Following steps similar to those used to prove Theorem \ref{thm:network_nodes_general} we get
\begin{equation}\label{eqn:LV_2}
V_e(e)\tilde F(t,e) \le 2(K_f-\sigma\lambda_2) V(e).
\end{equation}
The next step to estimate $LV(e)$ is that of computing $\frac{1}{2}tr\left\{ \tilde G(e)^TV_{ee}\tilde G(e)\right\}$. That is,
$$
\left(tr\left\{ (\sigma^\ast)^2 e^T((L^\ast)^T\otimes I_n)(L^\ast\otimes I_n)e\right\}\right)^{1/2} =\sigma^\ast\norm{(L^\ast\otimes I_n)e}_F = \sigma^\ast\abs{(L^\ast\otimes I_n)e}.
$$
Therefore, we have:
$$
\frac{1}{2}\left(tr\left\{ \tilde G(e)^TV_{ee}\tilde G(e)\right\}\right) = \frac{1}{2}(\sigma^\ast)^2\abs{(L^\ast\otimes I_n)e}^2 = \frac{1}{2}(\sigma^\ast)^2 e^T((L^\ast)^T\otimes I_n)(L^\ast\otimes I_n)e.
$$
We finally have:
\begin{equation}\label{eqn:VXX_2_2}
\frac{1}{2}tr\left\{ \tilde G(t,e)^TV_{ee}\tilde G(t,e)\right\} \le \frac{1}{2}(\sigma^\ast)^2 (\lambda_N^\ast)^2e^Te = (\sigma^\ast)^2(\lambda_N^\ast)^2 V(e).
\end{equation}
Combining (\ref{eqn:VXX_2_2}) and (\ref{eqn:LV_2}) we get:
\begin{equation}\label{eqn:c2}
LV(e) \le \left(2\left(K_f-\sigma\lambda_2\right) + (\sigma^\ast)^2(\lambda_N^\ast)^2\right)V(e).
\end{equation}
{\bf Estimate of $\abs{V_e(e)\tilde G(t,e)}^2$.} By means of Lemma \ref{lem:laplacian} we have:
$$
\abs{V_e(e)\tilde G(t,e)} = \sigma^\ast e^T(L^\ast\otimes I_n)e \ge \sigma^\ast\lambda_2e^Te.
$$
Hence
$$
\abs{V_e(e)\tilde G(t,e)}^2 \ge (\sigma^\ast)^2(\lambda_2^\ast)^2(e^Te)^2 = 4(\sigma^\ast)^2(\lambda_2^\ast)^2V(e)^2.
$$
Therefore, the network achieves complete stochastic synchronization if
$$
4(\sigma^\ast)^2(\lambda_2^\ast)^2> 2\left(2\left(K_f-\sigma\lambda_2\right) + (\sigma^\ast)^2(\lambda_N^\ast)^2\right),
$$
which is true by hypotheses.
\endproof
\subsubsection*{Example: field potentials revised}
\begin{figure}[tbh]
\begin{center}
\centering
\psfrag{x1}[c]{{$\gamma$}}
\psfrag{x2}[c]{{$t$}}
\psfrag{x3}[c]{{$t$}}
\psfrag{y1}[c]{{$R$}}
\psfrag{y2}[c]{{$v_i$'s}}
\psfrag{y3}[c]{{$v_i$'s}}
  \includegraphics[width=8cm]{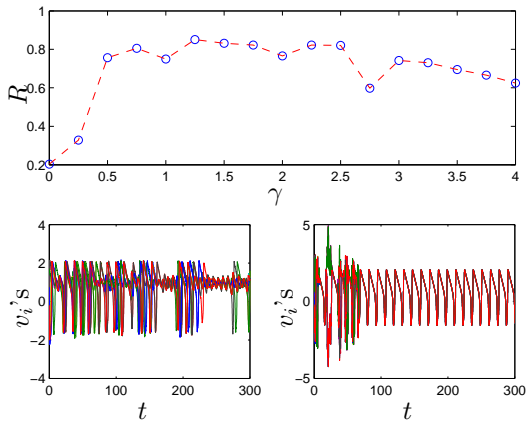}
  \caption{Simulation of (\ref{eqn:FN_full}) with the same parameters used in Section \ref{sec:com_sto_syn}. In the top panel, the behavior of $R$ is shown as a function of the noise intensity, $\gamma$. Such a figure shows a transition to synchronization caused by the increase in $\gamma$. This is confirmed by the bottom panels, where the time behaviors of the network state variables are shown for $\gamma = 0.25$ (left-bottom panel) and $\gamma = 3.75$ (right-bottom panel). The simulations were run using Matlab/Simulink, with solver method {\it ode23s} and relative tolerance {\it 1e-5}. Note that the decrease in $R$ for higher $\gamma$ is caused by the evaluation of the order parameter during the transient time.}
  \label{fig:order_param_complete}
  \end{center}
\end{figure}

In this Section we will apply Theorem \ref{thm:control} to show how noise diffusion can be used to synchronize a network. Specifically, we will again consider a network of FN oscillators coupled through noisy field potentials but this time coupling will be on both the nodes' state variables. In this case, the resulting stochastic differential equation is:
\begin{equation}\label{eqn:FN_full}
\begin{array}{*{20}l}
d v_i = \left[c\left(v_i+w_i -\frac{1}{3}v_i^3 +u(t)\right)\right]dt + \left[\gamma\left(\frac{1}{N}\sum_{j=1}^N(v_j) - v_i\right)\right]db\\
d w_i = \left[-\frac{1}{c}\left(v_i-a+bw_i\right)\right]dt+ \left[\gamma\left(\frac{1}{N}\sum_{j=1}^N(w_j) - w_i\right)\right]db, & i=1,\ldots,N.\\
\end{array}
\end{equation}
In terms of the notation introduced in Theorem \ref{thm:control} we have (please see the Appendix for the computations): (i) $\sigma = \lambda_2 = 0$; (ii) $K_f = \frac{-b+c^2 +\sqrt{1+b^2+2(b-1)c^2+2c^4}}{2c}$; (iii) $\sigma^\ast = \gamma$, $\lambda_N^\ast = N$ and $\lambda_2^\ast = N$. Therefore, following Theorem \ref{thm:control}, we have that network (\ref{eqn:FN_full}) will synchronize if $\frac{N^2\gamma^2}{2}>K_f$. That is, the network will synchronize if noise is strong enough and/or the number of nodes is sufficiently large. In Figure \ref{fig:order_param_complete} we simulated network (\ref{eqn:FN_full}) with $N = 10$.  In such a figure (top panel) the order parameter for network (\ref{eqn:FN_full}) is plotted as a function of the noise intensity $0\le \gamma \le 4$. Such a panel shows that the increase in $\gamma$ causes a transition from an unsynchronized state to synchronization. This is also confirmed by the bottom panels of Figure \ref{fig:order_param_complete} where the time behaviors of $v_i$'s are shown for $\gamma = 0.25$ and $\gamma = 3.75$.

\section{Application to collective decision making}\label{sec:applications}

We now further show the effectiveness of our results by considering the problem of understanding whether a collective decision can be made in group of agents. Over the past few years, the study of collective decision making processes (also known as collective cognition) has received significant attention, \cite{Con_Lis_09}. In particular, human performance in a two alternative decision making process is well modeled by the so-called Drift Diffusion Model (DDM). In its simplest instance, a DDM has the form (see e.g. \cite{Rat_McK_08}, \cite{Pet_Sri_Tay_Sur_Eck_Bul_15}  and references therein) $dx = \beta dt + \sigma dB$,
where $\beta \in\R$ is the drift rate, $\sigma >0$ is the diffusion rate and $x(t)$ is the aggregate evidence (i.e. opinion) at time $t$. Recently, in \cite{Sri_Leo_14} the speed-accuracy trade-off in collective decision making processes has been studied. In order to achieve their results, the authors used the coupled DDM, described by the  stochastic differential equation $dx_i = \left[\beta + \sum_{j\in\sN_i}(x_j-x_i)\right]dt + \sigma db_i
$, where: (i) $x_i$ is the aggregate evidence of the $i$-th agent in the network; (ii) $db_i$ models the (external) noise affecting the data collected by $i$-th agent; (iii) $\sigma$ models the strength of noise diffusion on each network node. Note that all the agents in the network have the same drift rate, $\beta$. 

In this Section we will consider a variation of the coupled DDM, where the noise diffusion at agent $i$ depends on how far the agent's opinion is with respect to the overall group agreement. We will also consider a state-dependent nonlinear drift for each agent rather than the constant drift rate $\beta$. The drift rate we will consider for the $i$-th agent is $x_i - x_i^3$ and models the fact that each of the agents in the network has the possibility to choose between two mutually excluding alternatives. This is  formalized with the following stochastic differential equation 
\begin{equation}\label{eqn:dec_making}
dx_i = \left[x_i - x_i^3 + \sum_{j\in\sN_i}(x_j-x_i)\right]dt + \left[x_i-\frac{1}{N}\sum_{j=1}^Nx_j\right]db_i.
\end{equation}
Equation (\ref{eqn:dec_making}) also arises in many situations from the IoT. For example, when, based on the local collection of data (observations), a group of interconnected objects with differing likelihoods is in charge of detecting the occurrence of a given event among a set of alternatives.

\begin{figure}[tbh]
\begin{center}
\centering 
\psfrag{x}[c]{{$t$}}
\psfrag{y}[c]{{$x_i$'s}}
\includegraphics[width=4cm]{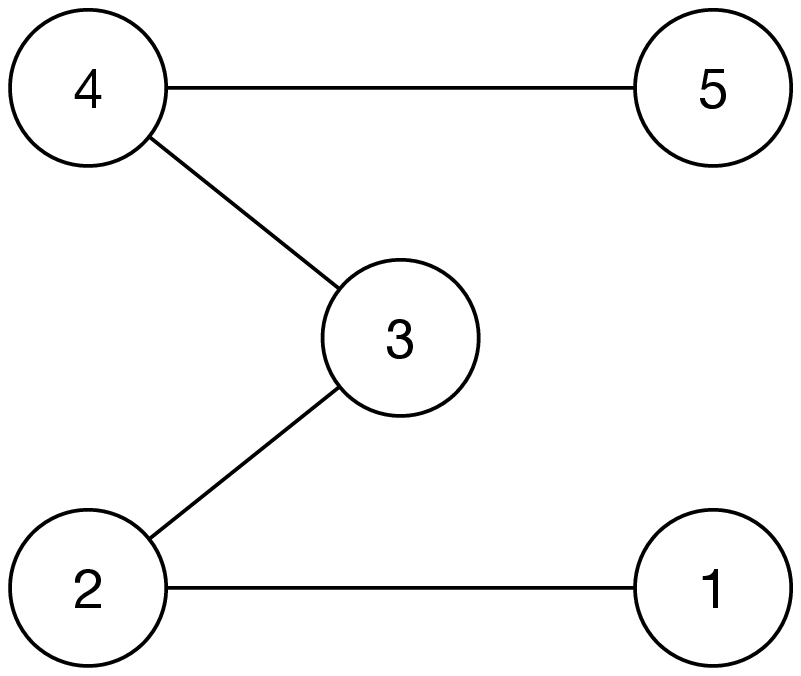}
\includegraphics[width=4cm]{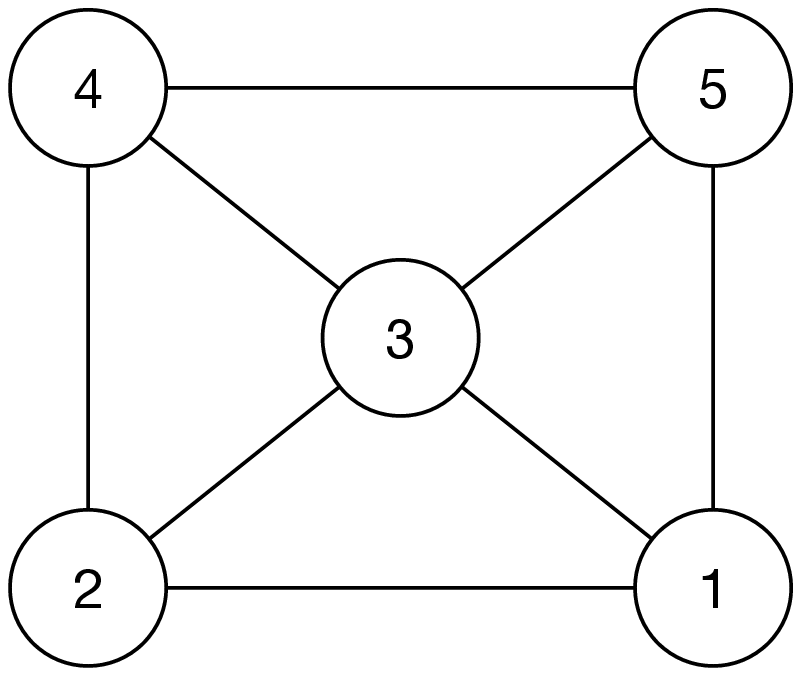}
\includegraphics[width=4cm]{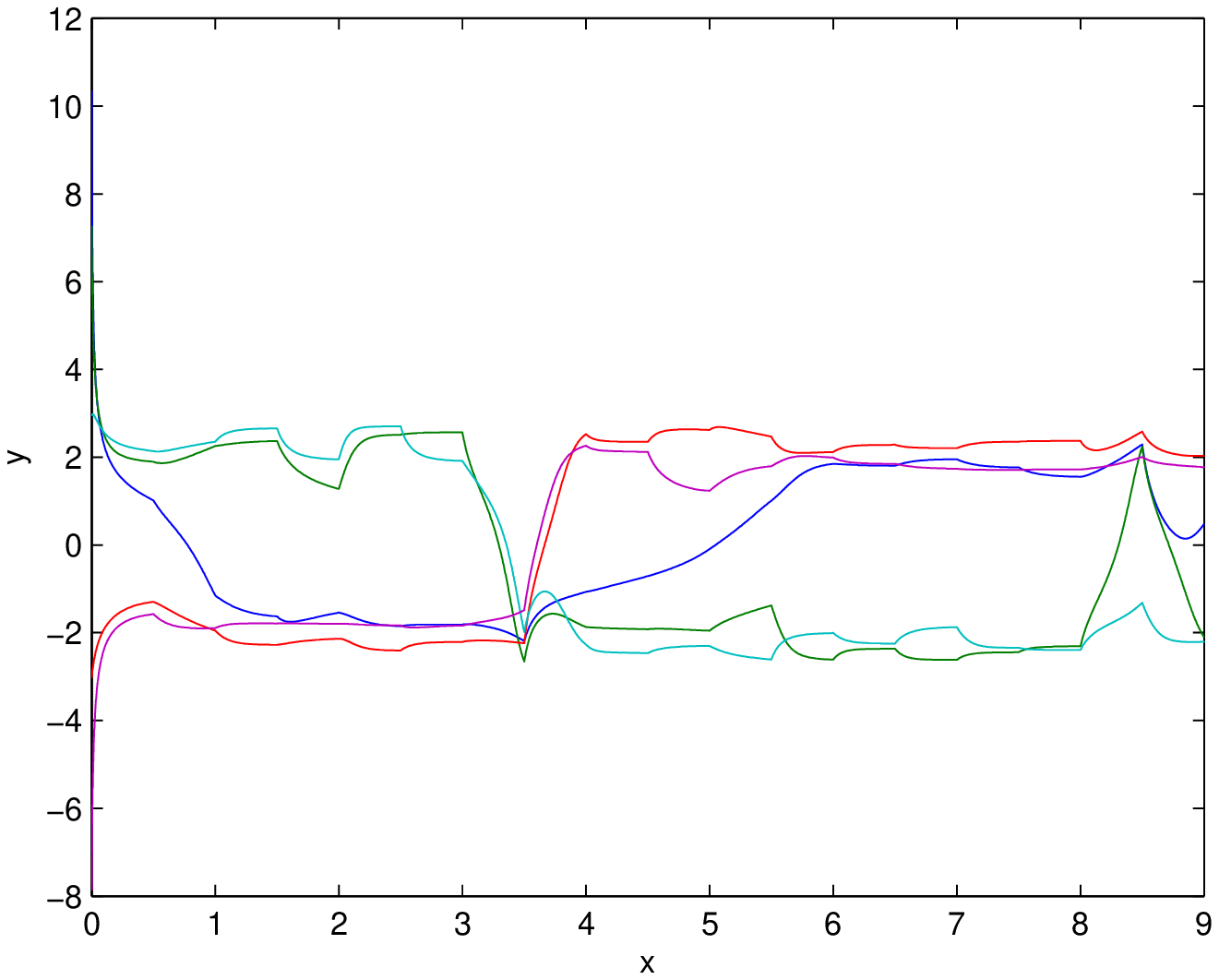}
\includegraphics[width=4cm]{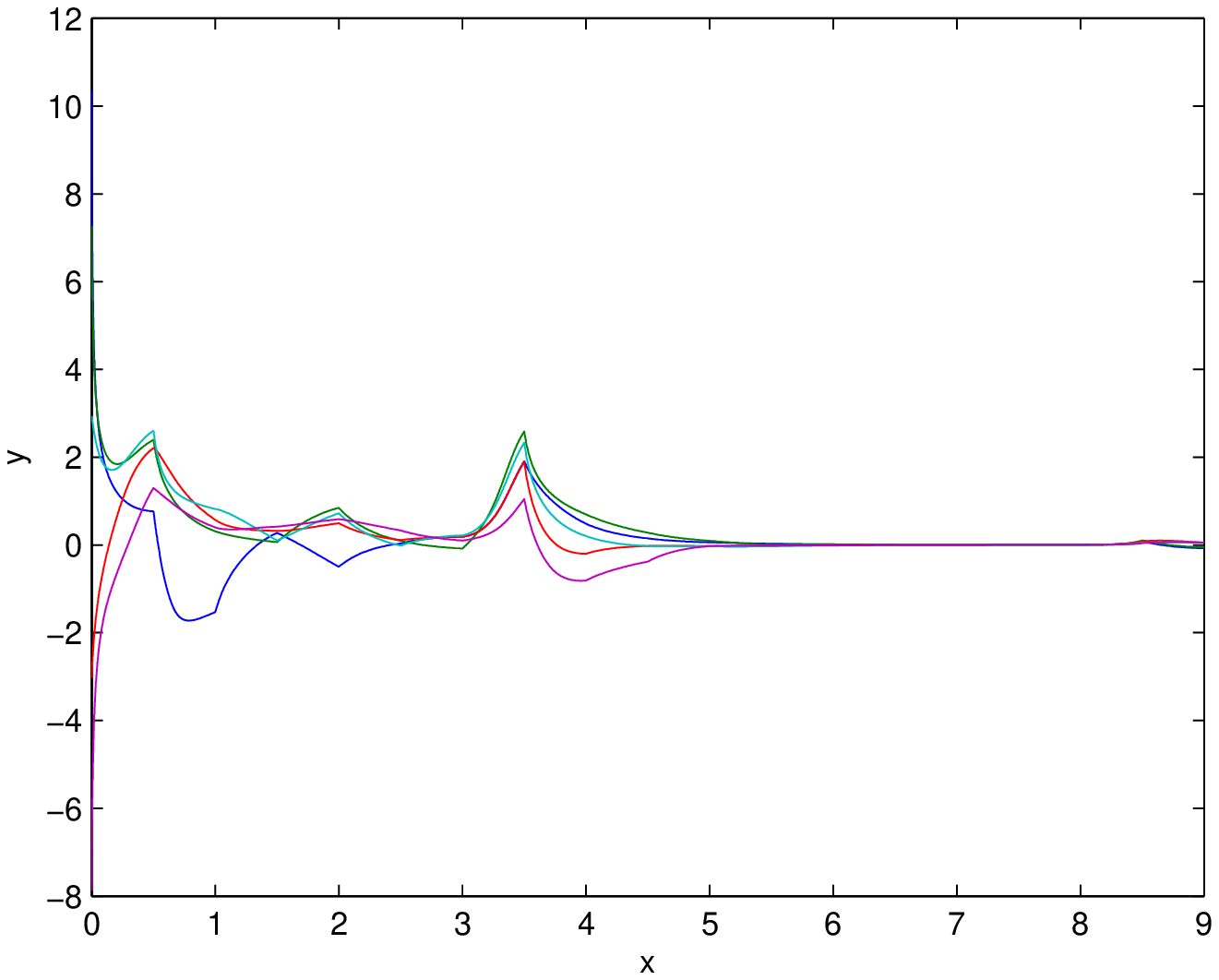}
  \caption{The two network topologies considered for the decision making process. Topology in the top-left panel has $\lambda_2 = 0.38$, while for the topology in the top-right panel $\lambda_2 = 3$. As shown in the time series at the bottom panels, collective decision is achieved only for the latter network. Initial conditions are chosen from a normal distribution of normal deviation equal to $10$ ($x_0 = [10.3469,   7.2689,-3.0344,2.9387,-7.8728]$).}
  \label{fig:net_topologies}
  \end{center}
\end{figure}

We will now use Theorem \ref{thm:network_nodes_general} to devise a condition ensuring that all network nodes achieve a common collective decision. As shown in the Appendix, application of Theorem \ref{thm:network_nodes_general} yields to the conclusion that (\ref{eqn:dec_making}) will achieve a collective decision if $\lambda_2 >0.5$. In order to confirm this theoretical prediction, we considered the two networks of Figure \ref{fig:net_topologies}. The topology in the top-left panel of Figure \ref{fig:net_topologies} violates the condition provided by Theorem \ref{thm:network_nodes_general} ($\lambda_2 = 0.38$). The bottom-left panel of the same figure shows that no collective decision is achieved by the network but rather network nodes take different decisions. The topology in the top-right panel of Figure \ref{fig:net_topologies} instead fulfills the condition provided by Theorem \ref{thm:network_nodes_general} ($\lambda_2 = 3$). In this case, as shown in the bottom-right panel, network nodes all converge towards an {\em agreed} behavior, indicating that a collective decision is achieved.

\section{Conclusions}
In this paper we presented new sufficient conditions ensuring that a network modeled via an Ito stochastic differential equation achieves complete stochastic synchronization. In Section \ref{sec:com_sto_syn} we first provided a condition on the coupling strength and network topology ensuring that the network synchronizes {\em despite} noise. Then, in Section \ref{sec:noise_synch}, we showed that, if properly designed, noise can be actually made useful for synchronization. Specifically, we showed that an additional layer can be designed where noise can propagate to {\em drive} the network towards a synchronous behavior. At the best of our knowledge, this approach to control synchronization is new. Our results have been strongly motivated by our current research in synchronization mechanisms of noisy networks and control applications in smart cities and IoT. In the former case, we are mainly interested in understanding whether noise propagation can destroy/enhance synchronization. In the latter case, we are mainly interested in controlling synchronization when the coupling between nodes needs to be low and when privacy between nodes needs to be preserved. We showed the effectiveness of our approach by using our results to analyze/control the onset of synchronization in networks of FN oscillators and to study collective decision processes. Based on the results of this paper, some of the future directions for our research will include: (i) extending our results to switching and directed topologies; (ii) use noise to solve distributed optimization problems; (iii) study social dynamics.

\appendix
\section*{Appendix}

\subsection*{Parameters for the example of Section \ref{sec:com_sto_syn}}
{\bf Computation of $K_f$.} The symmetric part of the Jacobian of the FN dynamics is:
$$
J_{sym}(v) = \left[ {\begin{array}{*{20}{c}}
   {c - {v^2}} & {{{\left( {c - \frac{1}{c}} \right)} \mathord{\left/
 {\vphantom {{\left( {c - \frac{1}{c}} \right)} 2}} \right.
 \kern-\nulldelimiterspace} 2}}  \\
   {{{\left( {c - \frac{1}{c}} \right)} \mathord{\left/
 {\vphantom {{\left( {c - \frac{1}{c}} \right)} 2}} \right.
 \kern-\nulldelimiterspace} 2}} & { - \frac{b}{c}}  \\
\end{array}} \right]
$$
From trivial matrix algebra, we have that $\lambda_{max}(J_{sym}(v))\le\lambda_{max}(J_{sym}(0))$, where
$\lambda_{max}(J_{sym}(0))=\frac{-b+c^2+\sqrt{1+b^2+2(b-1)c^2+2c^4}}{2c}$.
Since $u(t)$ is a smooth function, we have that there exists a $\tilde{\xi}\in[0,1]$ such that
$(x-y)^T[f(x,t)-f(y,t)]=(x-y)^T\left(\frac{\partial}{\partial x}{f\left(y+\tilde\xi(x-y),t\right)}\right)(x-y)$, and thus $K_f = \lambda_{max}(J_{sym}(0))$.

{\bf Computation of $K_g$ and $\bar K_g$.} Let $x_i:= (v_i, w_i)^T$ and $X:=[x_1^T,\ldots,x_N^T]^T$. Then, writing the dynamics of network (\ref{eqn:FN}) in compact form yields to $G(t,X) = G(X) = -(A\otimes \Gamma)X$, where: (i) $\Gamma$ is the $2\times 2$ matrix having $\Gamma_{11} = 1$ and all of its other elements equal to $0$; (ii)
 $A = I_N - 1_{N\times N}$. Let $e = X-S$, we have $G(X) - G(S) = -(A\otimes \Gamma)e$. Now, $K_g$ and $\bar K_g$ are obtained by noticing that the smallest eigenvalue corresponding to the direction transversal to the synchronization manifold is $0$, while the largest eigenvalue is $N$.

\subsection*{Parameters for the example of Section \ref{sec:noise_synch}}
Computation of $\lambda_2^\ast$ and $\lambda_N^\ast$ immediately follows from algebraic graph theory (see e.g. \cite{God_Roy_01}). Specifically, it suffices to note that the Laplacian of the noise diffusion layer has all the eigenvalues corresponding to the eigenvectors transverse to the synchronization manifold equal to $N$.

\subsection*{Parameters for collective decision making}
It is straightforward to check that $\sigma = 1$. In order to compute $K_f$, note that $(x-y)(f(x) - f(y))$ in this case amounts to $(x-y)((x-y)-(x^3-y^3))$ and therefore $(x-y)(f(x) - f(y)) \le (x-y)^2$. That is, $K_f = 1$. Note that writing (\ref{eqn:dec_making}) in compact form yields to the diagonal matrix $G(X)$ having $g_{ii} = x_i - \frac{1}{N} \sum_{j=1}^Nx_j = x_i -s = e_i$. Therefore, $\norm{G(X)}_F = \sqrt{\sum e_i^2}$, thus implying that $K_g = 1$. Moreover, $e^TG(X) = \sum e_i^2$ and therefore $\bar K_g = 1$.

\end{document}